# A Semantic Situation without Syntax (Non- axiomatizibility of Theories)

Farzad Didehvar
didehvar@aut.ac.ir

Department of Mathematics &Computer Science,
Amir Kabir University, Tehran, Iran
Tel: 98-21-66952518

**Abstract:** *Here, by introducing a version of "Unexpected hanging paradox" first we try to open a new way and a new explanation for paradoxes, similar to liar paradox. Also, we will show that we have a semantic situation which no syntactical logical system could support it. Finally, we propose a claim in the subject of axiomatizibility. Based on this claim, having an axiomatic system for Computability Theory is not possible. In fact, the same argument shows that many other theories are non-axiomatizable. (Dare to say: General Theories of Physics and Mathematics).*

**Introduction:** *The surprise exam paradox (Unexpected hanging paradox) is a well-known paradox that too many works have been done based on that or about it.*

*Simple formulation of this paradox is:*

*A judge condemned a prisoner. The verdict is: you will be executed in the next week, in a day that you did not know about it in previous days, so it will be a surprise for you.*

*The prisoner starts to argue as following:*

1. *I will not be executed in Friday, since in this case, I will not be executed in the days before so in the Thursday I will understand that I will be executed in Friday. A contradiction.*
2. *I will not be executed in Thursday, since in this case; I will not be executed in the days before by 1 I know that I will not be executed in Friday, so in the Wednesday I will understand that I will be executed in Thursday. A contradiction.*
3. *Analogously I will not be executed in Wednesday, Tuesday, Monday and Sunday.*

*He concludes that he will not be executed in next week. But the guard executed him in Tuesday. The guard explanation was, as he claimed he was not expected to be executed in the*

*Next week, so we did not contradict the verdict.*

*This paradox is known as "Surprise test paradox" too. For some technical reasons we do not apply this name here and we use the more origin name. "The bottle Imp" and "Crocodile dilemma" has some similarity to this paradox.*

*Historically, O'Conner article [14] was the first article which expose this paradox academically.*

*Some people think that basically, it is not a paradox . Quine work [15] was in this line. Some think the paradox has been resolved by Quine And simply some aspects of it remains to be discussed [8].*

*In fact, there is no consensus on its nature so we have no final resolution yet, besides all we have different various reformulation of this paradox.*

*There are different approaches to solve this problem:*

*1. <span style="color:red">Logical approach :</span> But some stress more on its logical aspects and they try to attack this problem by a logical approach. As an example some thinks to apply and to assume self contradictory and self referencing the nature of this problem.* <span style="color:green">*(referencing some works)*</span>

*Pararel to this approch, this paradox is applied to give a solution to the first Gödel Theorem and the Second Gödel Theorem in [6] inspired by surprise examination argument and applying some other Theorems the authors give a new proof for second completeness theorem. Before this work chaitin by applying Berrys paradox has done a similar work [4].*

*It is considered as an attempt to show the self referential features of this paradox.*

*There is an similar elder article, by [3].*

*The explanation of [16] is in this line.*

*2. Epistemic approach: Some think about this paradox as a paradox of knowledge (Epistemic Paradox),In some approaches it goes to ability or inability to know some statement[9] .In recent years we have some attempts to solve this problem in dynamic epistemic logic.[12], [13]*

*There are attempts yet to solve this problem once and forever.*

*While many solutions have been offered to the surprise exam paradox over the past sixty*

*Years it continues to perplex and torment. This paper attempts to end the intrigue once and for all.[11]*

*In this article by adding together three previous proof the author tries to give a new proof[unprojectible announcement argument],[wright &Sudbury argument],[epistemic blind spot argument].Too many philosophical article is published related to this subject.*

*Here, throughout this article, our goal is totally different and the article is not in this line at all .We apply a version of this paradox to show that there is a semantic situation that no syntax supports it.*

*In [5] by applying a version of "Unexpected hanging paradox", we tried to show that there is a proof which doesn't show the truth. Here, we try to modify and formalize this problem and by applying a version of this paradox to show that there is a semantic situation no syntax supports It.*

*In conclusion part, the above claim shows us a general claim about unaxiomatizability for so many Theories, as it could be used somehow as a criteria of axiomatizability of theories.*

In [5] by applying a version of "Unexpected hanging paradox", we tried to show that there is a proof which doesn't show the truth. Here, we try to modify and formalize this problem *and by applying a version of this paradox to show that there is a semantic situation no syntax supports It. .* <span style="color:red">Some similar attempts could be found in [1].</span>

*In conclusion part, the above claim shows us a general claim about unaxiomatizability for so many Theories, as it could be used somehow as a criteria of axiomatizability of theories.*

First, we explain this version of scenario of paradox.
A Computer Scientist invents a Computer which argues logically, and also any information and conclusion of this Computer is announced to public openly.

The People of the city do not like this Computer and they appoint a jury which decides about the destiny of this Computer.

The jury decision is declared to Computer Scientist and Computer as follows:

The Computer will be destroyed in the next week on a day that it wasn't concluded by the computer itself before that day. Definitely we know everything that the computer knows or everything it will conclude. (The computer is invented in a way that he/it declares anything that he concludes; we can imagine that he gives all of his conclusions in a printing form).

After announcing it, the computer starts to argue and it finds the following arguments (the argument of 'hanging paradox'). By the way, because of Intelligence of the computer, all the people know destroying the computer is equal to "executing the Computer".

The Demonstration of a Computer:

I will not be executed on Saturday, since if I was executed on Saturday I would be alive till Friday,

and on Friday I conclude that I will be executed on Saturday ( since based on what jury said I will be executed in the next week). So in the case if I remain alive (!) till Friday, I will conclude on that day, I will not be executed on Saturday! (considering the proposed claim of the jury).This argument exists in my data base right now (after concluding by me, right now).In a more exact form, the Computer argues as follows: I have a timer, by using my storage and database I conclude that I will conclude on Friday "I will be executed on the next day, so based on the claims of juries, I will not be executed on Saturday".

2. I will not be executed on Friday, since if I was executed on Friday I would

remain alive till Thursday, and on Thursday and our claim in point 1 and based on what jury said("The Computer will be executed in the next week") I conclude that I will be executed on Friday. So by the proposed claim of the jury I will not be executed on Friday!(This argument right now exists in my data base after concluding by me, so in the case if I remain alive(!) till Thursday, I will conclude on that day, that I will not be executed on Friday). In a more exact form, the Computer argues as follows: I have a timer, by using my storage and database I conclude that I will conclude "on Thursday I will be executed on the next day, so by the proposed claim of juries, I will not be executed on Friday".

…

6. I will not be executed in Monday, since if I will be executed in Monday I will be alive till Sunday,

And on Sunday and the claims in point numbers 1&2&…&5, I conclude that I will be executed on Monday, since by what jury said I will be executed in the next week. So in case if I remain alive (!) till Sunday, I will conclude on such a day that I will not be executed on Monday! (On the basis of proposed claims of the jury), this argument exists right now in my data base after concluding by me. In a more exact form, the Computer argues as follows: I have a timer, based on my storage and database I conclude that I will conclude on Sunday "I will be executed on the next day, so by the claim of juries, I will not be executed on Monday!".

7. I will not be executed on Sunday, since based on points 1-6 I will not be executed on Saturday, Friday, … Monday .

So I should be executed on Sunday! My conclusion (to be executed on Sunday) concludes that I will not be executed on Sunday!

So, I will not be destroyed in the next week.

The computer didn't conclude more, except his every day conclusion he concluded he will not be executed on the next day based on the similar above argument.

The computer was destroyed on Tuesday (has was executed), and the Computer Scientist complained about this injustice, and he transferred this argument of Computer to journals and the court. In this message the computer scientist mentions that the computer presented a logical argument, based on that he will not be executed in the next week, nevertheless he was executed).

The court said:

The Computer proved that he will not be executed. This is a true proof. We executed him on Tuesday, and as he claimed he didn't conclude that he will be executed on Tuesday. On the contrary, he announced that "he will not be executed on that day". In other words, he proved in his message that by accepting what the judge said as a true claim, he would not be executed in that week. More formally, here we state as follows:

P: we consider what the judge said as a true claim

Q: He would not be executed in this week

"His proof (p|---q) is a true proof, but it doesn't show the Truth". (*)

In this paper, we defend the above claim of the juries (*) and we try to show how we could develop this idea.

It is noteworthy to say that, we don't claim that the above result does not propose a solution for unexpected hanging paradox in all its versions. It is simply considered as the only way to explain this version of "Unexpected hanging paradox". Later on, we will know the above result as a possible explanation for the other versions of this paradox and some other paradoxes like liar paradox. In [2], [4] we can find a more conclusive explanation of Liar paradox.

To formalize the above proof, we call our Computer "A" and whenever we say "A concludes", in our paradox, we replace it by: A:$[\varphi]$. Also, A:$_i[\varphi]$ stands for "A concludes $\varphi$ in the i$^{th}$ day".

In the following formalism, $\varphi(i)$ stands for A will be executed in the i$^{th}$ day. Now in any formalization of the problem and proof we will have the following assertions:

1. (A:$_i[\varphi(i+1)]) \to \sim\varphi[i+1]$ (If A utters in the i-th day that he will be executed in i+1 th day, he will not be executed in i+1 th day, i=0, 1,...,6).
2. $\vee_{i=1}^{i=7} \varphi(i)$
3. A:$_i[\sim\varphi(i+1)] \to (A:_i [\varphi(i+1)])$ ( If A says that I will not be executed in tomorrow, he will not say that, he will be executed( Since A argues logically)).
4. $\varphi(k+1) \to \{\wedge_{i=1}^{k} \sim\varphi(i)\}$(If he is executed in a day, we conclude the days before he wasn't executed).
5. $\varphi(3)$.

It is notable that any syntactical system for the above paradox includes the above statements, directly or as a conclusion.

**Proof:**

1. $\sim\varphi(7)$.
   If $\varphi(7)$ by 4 we have $\wedge_{i=1}^{6} \sim\varphi(i)$. By 2 we have A:$_6[\sim\varphi(7)]$ . By 3 we have (A:$_6 [\varphi(7)]$). By 1 we have
   $\sim\varphi[7]$.
2. $\sim\varphi(6)$. Similar to above proof.
   
   ......
   
   ......

   1. $\sim\varphi(1)$. In a similar way.
   2. So we have a contradiction here, since $\varphi(3)$.

So, we have a model which its associated formal system is contradictory!

It is notable that any syntactical system that we attribute to the above paradox contains the above statements, whether directly or as a conclusion.

In other words, the proof shows that this formalism and any other formalism are essentially contradictory but we have a semantic situation for this contradictory formalism.

**Remark:** In any complete formalization of the above paradox we need to extend our language in a way the language includes some sentences in the form of:

$[A_{:i} \ [A_{:j} \ (\varphi)]]$,

In order to write some statements like the following statement in the demonstration of Computer, I conclude that I will conclude on Friday "I will be executed on the next day, so by the claim of juries, I will not be executed on Saturday".

Actually we have the following axioms:

$[A_{:0} \ [A_{:j} \ (\varphi(j+1))]] \to [A_{:j} \ (\varphi(j+1))]$ J=0,…, 6.

However, for simply showing the contradiction existed in formalizing this semantic situation, it is not essential to face this extension of the language, and the above axioms.

**An explanation and conclusion:**
In the above system, there is a contradiction. In brief, the judges claim that the Computer will be executed next week, but the Computer proves that he will not. So this system is a contradictory system, but at the same time we have a semantic situation for this system. So we have a contradictory system which has a semantic situation.
 In other word, we have semantic situation which no consistent syntactical system supports them. (Here, we have such an example).
As a result, our proof in above does not show any truth. So, there are some incorrigible flaw in modeling and formalizing the proofs. In other word, formal systems are not able to support such semantic situations.

Clearly, this opens a new way to explain some paradoxes similar to liar paradox, as follows:

In such paradoxes the proofs don't show the truth, since there is no consistent syntactical system to support the related semantic situation.

This would be considered as the central result and theme of this paper. As the last word about the subject of Computability Theory (and Theories in general), we propose the following conclusion.

- Conclusion: *In the above formalism, we considered A as a Turing machine. By a slight modification in the formalism, we are able to replace A by [A] (code of A), in a similar way $\varphi$ by $[\varphi]$. So the above problem has a formalism in the scope of Computability Theory. If one of the real goals of Computability Theory is to explain the real situations related to the Computations, the Computability Theory should explain the above semantic situation. Since any axiomatization of the semantic situation falls in contradictory, we know that no axiomatization is able to deal with the semantic situation. So Computability Theory is not axiomatizable. Actually, if Computability Theory does not explain the above situation; it would be a so naïve Theory and will contradict our expectations.*

*We are able to generalize the above assertion to any other theory in which we expect to explain the formalized form of the scenario.*

*It is rational to say:*

*"In the Theories of Physics and Mathematics, we have such an expectation that these theories will be able to explain the scenario, but in Euclidean Geometry and Theory of probability, we do not expect these theories to explain the situation of our scenario, but they are axiomatizable theories. Regarding practical objectives, axiomatazibility of some other theories are recognizable by means of our scenario.*

In other words, practically it is useful to ask the following question when we face a Theory:

Do we expect the proposed Theory explains the situation of our scenario?

If we expect so, our Theory will not be axiomatizable.